\documentclass[11pt]{amsart}

\usepackage{amsthm, amsfonts, amssymb, amscd}
\usepackage[
colorlinks]{hyperref}
\usepackage{tikz-cd}
\usepackage{geometry}
\usepackage{marginnote}

\theoremstyle{definition}
\newtheorem{ntn}{Notation}[section]

\theoremstyle{plain}

\newtheorem{prp}[ntn]{Proposition}
\newtheorem{thm}[ntn]{Theorem}
\newtheorem{cor}[ntn]{Corollary}

\theoremstyle{definition}

\newtheorem{rem}[ntn]{Remark}
\newtheorem{exa}[ntn]{Example}

\numberwithin{equation}{section}

\newcommand{\N}{\mathbb{N}}
\newcommand{\z}{\mathbb{Z}}
\newcommand{\q}{\mathbb{Q}}

\newcommand{\C}{\mathbb{C}}

\newcommand{\F}{\mathbb{F}}

\newcommand{\OO}{\mathcal{O}}

\newcommand{\mmm}{\mathfrak{m}}

\renewcommand{\aa}{{A^\times}}

\newcommand{\pth}{{\Big[\frac{1}{p}\Big]}}

\newcommand{\lan}{\langle}
\newcommand{\ran}{\rangle}
\newcommand{\se}{\subseteq}
\newcommand{\arr}{\rightarrow}
\newcommand{\larr}{\longrightarrow}

\newcommand{\two}{\twoheadrightarrow}

\newcommand{\ep}{\varepsilon}

\newcommand{\GL}{\mathit{{\rm GL}}}
\newcommand{\Ee}{\mathit{{\rm E}}}
\newcommand{\SL}{\mathit{{\rm SL}}}

\newcommand{\St}{\mathit{{\rm St}}}
\newcommand{\PSL}{\mathit{{\rm PSL}}}

\renewcommand{\char}{{\rm char}}
\newcommand{\diag}{{\rm diag}}

\newcommand{\Vv}{{\rm V}}
\newcommand{\Kk}{{\rm K}}
\newcommand{\Uu}{{\rm U}}
\newcommand{\Cc}{{\rm C}}
\newcommand{\GE}{{\rm GE}}

\newcommand{\ab}{{\rm ab}}

\newcommand {\mtxx}[4]
{\left(\!\!
	\begin{array}{cc}
		\!\!#1 & \!\!#2   \\
		\!\!#3 & \!\!#4
	\end{array}\!\!
	\right)
}

\newtheoremstyle{athm}
{}
{}
{\itshape}
{}
{\scshape}
{}
{.5em}
{\thmnote{#3}}
\theoremstyle{athm}

\begin{document}
	\title{The abelianization of the elementary group of rank two}
	\author{Behrooz Mirzaii}
	\author{Elvis Torres P\'erez}
	
	\begin{abstract}
		For an arbitrary ring $A$, we study the abelianization of the elementary group $\Ee_2(A)$.
		In particular, we show that for a commutative ring $A$ there exists an exact sequence
		\[
		\Kk_2(2,A)/\Cc(2,A) \arr A/M \arr \Ee_2(A)^\ab \arr 1,
		\]
		where $\Cc(2,A)$ is the central subgroup of the Steinberg group $\St(2,A)$ generated by the Steinberg symbols and
		$M$ is the additive subgroup of $A$ generated by $x(a^2-1)$ and $3(b+1)(c+1)$, with $x\in A, a,b,c \in \aa$.
	\end{abstract}
	\maketitle
	
	Let $A$ be an associative ring with 1. Let $\Ee_2(A)$ be the elementary subgroup of the general linear 
	group $\GL_2(A)$. In the current paper, we study the abelianization  of $\Ee_2(A)$, i.e. the group
	\[
	\Ee_2(A)^\ab:=\Ee_2(A)/[\Ee_2(A),\Ee_2(A)].
	\]
	This group has been studied in the literature before. The first important result in this direction was proved by 
	P.M. Cohn. In his seminal paper \cite[Theorem 9.3]{cohn1966}, Cohn showed that if $A$ is quasi-free for 
	$\GE_2$ and $\aa$ is abelian, then
	\[
	\Ee_2(A)^\ab\simeq A/M,
	\]
	where $M$ is the additive subgroup of $A$ generated by $axa-x$ and $3(b+1)(c+1)$ with $x\in A$ and $a,b,c\in \aa$
	(see \cite[p. 10]{cohn1966} for a definition of quasi-free for $\GE_2$). Later Cohn generalised this result to rings 
	that are universal for $\GE_2$ \cite[Theorem 2]{cohn1968}. In particular he showed that for any 
	ring $A$, there is a homomorphism $A/M \arr \Ee_2(A)^\ab$.
	
	As the main result of this paper we show that the homomorphism $A/M \arr \Ee_2(A)^\ab$ is a part of an exact 
	sequence with interesting terms. In particular, we show that for any commutative ring $A$, there is
	an exact sequence of $\z[\aa/(\aa)^2]$-modules
	\[
	H_2(\Ee_2(A),\z) \arr \frac{\Kk_2(2,A)}{[\Kk_2(2,A), \St(2,A)]\Cc(2,A)} \arr A/M \arr \Ee_2(A)^\ab \arr 0.
	\]
	We refer the reader to Theorem \ref{H11} for a general statement over an arbitrary ring. For the definition of
	$\St(2,A)$, $\Kk_2(2,A)$ and $\Cc(2,A)$ see Section \ref{CA} below.
	
	Many of known results about the structure of $\Ee_2(A)^\ab$, that we could find, follow from 
	this exact sequence and  we generalize most of them (see for example \cite[Theorem]{menal1979}, 
	\cite[Theorem 9.3]{cohn1966}, \cite[Theorem 2]{cohn1968} \cite[Corollary 4.4]{stein1971}). 
	For example it follows immediately that if $A$ is universal for $\GE_2$, then $\Ee_2(A)^\ab\simeq A/M$. 
	Moreover, we show that for any square free integer $m$, we have
	\[
	\begin{array}{c}
		\Ee_2(\z[\frac{1}{m}])^\ab\simeq 
	\end{array}
	\begin{cases}
		0   & \text{if $2|m,\ 3|m$}  \\
		\z/3   & \text{if $2|m,\ 3\nmid m$}\\
		\z/4   & \text{if $2\nmid m,\ 3|m$}\\
		\z/12 &  \text{if $2\nmid  m,\ 3\nmid m$.}
	\end{cases}
	\]
	Some parts of the above isomorphism were known. But during the preparation of this article we could 
	not find this general result in the literature. After that this paper appeared on arXiv, Nyberg-Brodda
	informed us that he also proved the above isomorphism with a different method. His proof now 
	can be found in \cite{nyberg2024}.
	
	Finally, we study $\Ee_2(A)^\ab$ over the ring of algebraic integers of a quadratic filed $\q(\sqrt{d})$.
	When $d<0$, this group was calculated by Cohn in \cite{cohn1966}, \cite{cohn1968}. If $d>0$, then 
	this group is finite and we give an estimate of its structure.
	
	\section{Elementary groups of rank 2 and rings universal for \texorpdfstring{$\GE_2$}{Lg}}\label{CA}
	
	Let $A$ be a ring (associative with $1$). Let $\Ee_2(A)$ be the subgroup of $\GL_2(A)$ generated by the 
	elementary matrices
	$E_{12}(a):=\begin{pmatrix}
		1 & a\\
		0 & 1
	\end{pmatrix}$ and $
	E_{21}(a):=\begin{pmatrix}
		1 & 0\\
		a & 1
	\end{pmatrix}$, $a\in A$.
	The group $\Ee_2(A)$ is generated by the matrices
	\[
	E(a):={\mtxx a 1 {-1} 0}, \ \ \ \ a \in  A.
	\] 
	In fact
	\[
	E_{12}(a)=E(-a)E(0)^{-1}, \ \ \ \ \ E_{21}(a)=E(0)^{-1} E(a), 
	\]
	\[
	E(0)=E_{12}(1)E_{21}(-1)E_{12}(1).
	\]
	Let ${\rm D}_2(A)$ be the subgroup of $\GL_2(A)$ generated by diagonal matrices and 
	let $\GE_2(A)$ be the subgroup of $\GL_2(A)$  generated by ${\rm D}_2(A)$ and $\Ee_2(A)$. 
	For any $a\in \aa$, let 
	\[
	D(a):=\begin{pmatrix}
		a & 0\\
		0 & a^{-1}
	\end{pmatrix}\in {\rm D}_2(A).
	\]
	Since 
	\[
	D(-a)=E(a)E(a^{-1})E(a),
	\]
	the element $D(a)$ belongs to $\Ee_2(A)$. It is straightforward to check that
	\begin{itemize}
		\item[(1)] $E(x)E(0)E(y)=D(-1)E(x+y)$,
		\item[(2)] $E(x)D(a)=D(a^{-1})E(axa)$,
		\item[(3)] $\prod_{i=1}^l\! D(a_ib_i)D(a_i^{-1})D(b_i^{-1})\!=\!1$ 
		provided $\prod_{i=1}^l [a_i,b_i]\!=\!1$,
	\end{itemize}
	where $x,y\in A$ and $a, a_i, b_i\in \aa$. If $\aa$ is abelian, then (3) just is
	
	\begin{itemize}
		\item[(3$'$)] $D(ab)D(a^{-1})D(b^{-1})=1$ for any $a,b\in \aa$.
	\end{itemize}
	Note that $\Ee_2(A)$ is normal in $\GE_2(A)$ \cite[Theorem 2.2]{cohn1966}.
	
	Let $\Cc(A)$ be the group generated by symbols $\ep(a)$, $a\in A$, subject to the relations
	\begin{itemize}
		\item[(i)] $\ep(x)\ep(0)\ep(y)=h(-1)\ep(x+y)$ for any $x,y \in A$,
		\item[(ii)] $\ep(x)h(a)=h(a^{-1})\ep(axa)$, for any $x\in A$ and $a\in \aa$,
		\item[(iii)] $\prod_{i=1}^l\! h(a_ib_i)h(a_i^{-1})h(b_i^{-1})\!=\!1$ for any $a_i, b_i\!\in\!\aa$ 
		provided $\prod_{i=1}^l [a_i,b_i]\!=\!1$,\\
	\end{itemize}
	where 
	\[
	h(a):=\ep(-a)\ep(-a^{-1})\ep(-a).
	\]
	Note that by (iii), $h(1)=1$ and $h(-1)^2=1$. Moreover $\ep(-1)^3=h(1)=1$ and $\ep(1)^3=h(-1)$. 
	There is a natural surjective map 
	\[
	\Cc(A) \arr \Ee_2(A),  \ \ \ \ \ep(x)\mapsto E(x).
	\]
	We denote the kernel of this map by $\Uu(A)$. Thus we have the extension
	\[
	1 \arr \Uu(A) \arr \Cc(A) \arr \Ee_2(A) \arr 1.
	\]
	A ring $A$ is called {\it universal for} $\GE_2$ if $\Uu(A)=1$, i.e. the relations (i), (ii) and (iii) 
	form a complete set of defining relations for $\Ee_2(A)$. 
	
	\begin{exa}\label{example2}
		(i) Any local ring is universal for $\GE_2$ \cite[Theorem 4.1]{cohn1966}.
		\par (ii) Let $A$ be semilocal. Then $A$ is universal for $\GE_2$ if and only if none of the rings 
		$\z/2 \times \z/2$, $\z/6$ and ${\rm M}_2(\z/2)$ is a direct factor of $A/J(A)$, where $J(A)$ is the 
		Jacobson radical of $A$ \cite[Theorem 2.14]{menal1979}.
		\par (iii) If $A$ is quasi-free for $\GE_2$, then it is universal for $\GE_2$ (see \cite[p. 10]{cohn1966}).
		\par (iv) Any discretely normed ring is quasi-free for $\GE_2$ and hence is universal for $\GE_2$ 
		\cite[Theorem 5.2]{cohn1966} (see \cite[\S 5]{cohn1966} for the definition of a discretely normed ring). 
		\par (v) Let $\OO_d$ be the ring of algebraic integers of an imaginary quadratic field $\q(\sqrt{d})$ 
		($d<0$). Let $A$ be a subring of $\OO_d$. Then $A$ is discretely normed except for the rings $\OO_{-1}$, 
		$\OO_{-2}$, $\OO_{-3}$, $\OO_{-7}$, $\OO_{-11}$ and $\z[\sqrt{-3}]$ (see \cite[Propositions 1, 2]{dennis1975} 
		and \cite[\S 6]{cohn1966}). 
		\par (vi) Discretely ordered rings are universal for $\GE_2$ \cite[Theorem 8.2]{cohn1966}
		(see \cite[\S 8]{cohn1966} for the definition of a discretely ordered ring). If $A$ is discretely ordered, 
		then $A[X]$ is discretely ordered. The most obvious example of a discretely ordered ring is the ring $\z$. 
		Therefore $\z[X_1, \dots, X_n]$ is universal for $\GE_2$.
	\end{exa}
	
	
	\section{Rank one steinberg group} 
	
	The elementary matrices $E_{ij}(x)$, $i,j\in \{1,2\}$ with $i\neq j$, satisfy the following relations
	\begin{itemize}
		\item[(a)] $E_{ij}(r)E_{ij}(s)=E_{ij}(r+s)$ for any $r,s \in A$,
		\item[(b)] $W_{ij}(u)E_{ji}(r)W_{ij}(u)^{-1}=E_{ij}(-uru)$, for any $u\in \aa$ and $r\in A$,
	\end{itemize}
	where $W_{ij}(u):=E_{ij}(u)E_{ji}(-u^{-1})E_{ij}(u)$. Observe that 
	\[
	W_{12}(u)={\mtxx{0}{u}{-u^{-1}}{0}}, \ \ \ \ \ W_{21}(u)={\mtxx{0}{-u^{-1}}{u}{0}}.
	\]
	Set 
	\[
	H_{12}(u):=W_{12}(u)W_{12}(-1)=D(u), \ \ \ \ \ H_{21}(u):=W_{21}(u)W_{21}(-1)=D(u)^{-1}.
	\]
	
	The Steinberg group $\St(2,A)$ is the group with generators $x_{12}(r)$ and $x_{21}(s)$, $r,s\in A$, 
	subject to the Steinberg relations
	\begin{itemize}
		\item[($\alpha$)] $x_{ij}(r)x_{ij}(s)=x_{ij}(r+s)$ for any $r,s \in A$,
		\item[($\beta$)] $w_{ij}(u)x_{ji}(r)w_{ij}(u)^{-1}=x_{ij}(-uru)$, for any $u\in \aa$ and $r\in A$,
	\end{itemize}
	where 
	\[
	w_{ij}(u):=x_{ij}(u)x_{ji}(-u^{-1})x_{ij}(u).
	\]
	
	The natural map
	\[
	\Theta: \St(2,A) \arr \Ee_2(A), \ \ \ \ \ x_{ij}(r) \mapsto E_{ij}(r)
	\]
	is a well defined homomorphism. The kernel of this map is denoted by $\Kk_2(2,A)$ and is called the 
	unstable $\Kk_2$-group of degree $2$. For $u\in \aa$, let
	\[
	h_{ij}(u):=w_{ij}(u)w_{ij}(-1).
	\]
	It is not difficult to see that $h_{ij}(u)^{-1}=h_{ji}(u)$ \cite[Corollary A.5]{hutchinson2022}. 
	
	Let $u,v\in \aa$ commute and let
	\[
	\{u,v\}_{ij}:=h_{ij}(uv)h_{ij}(u)^{-1}h_{ij}(v)^{-1}.
	\]
	Since $\Theta(h_{12}(u))=H_{12}(u)=D(u)$ and $\Theta(h_{21}(u))=H_{21}(u)=D(u)^{-1}$, we have
	\[
	\Theta(\{u,v\}_{12})=D(uv)D(u)^{-1}D(v)^{-1}=0, \ \ \ \ \Theta(\{u,v\}_{21})=D(uv)^{-1}D(u)D(v)=0.
	\]
	Therefore, $\{u,v\}_{i,j}\in \Kk_2(2,A)$.
	
	The element $\{u,v\}_{ij}$ lies in the center of $\St(2, A)$  \cite[\S 9]{D-S1973}. It is straightforward to check that
	$\{u,v\}_{ji}=\{v,u\}_{ij}^{-1}$. We call $\{u,v\}_{ij}$ a {\it Steinberg symbol} and set
	\[
	\{v,u\}:=\{v,u\}_{12}=h_{12}(uv)h_{12}^{-1}(u)h_{12}(v)^{-1}.
	\]
	
	For a commutative ring $A$, let $\Cc(2,A)$ be the subgroup of $\Kk_2(2, A)$ generated by the Steinberg 
	symbols $\{u,v\}$, $u,v\in \aa$. Then $\Cc(2, A)$ is a central subgroup of $\Kk_2(2, A)$.
	
	For a ring $A$ let $\Vv_2(A)$ be the subgroup of $\aa$ generated by all elements $a\in \aa$ such that 
	$\diag(a,1)$ is in $\Ee_2(A)$. The following theorem is due to Dennis \cite[\S 9 (d), p. 251]{D-S1973}.
	
	\begin{thm}[Dennis]
		A ring $A$ is universal for $\GE_2$ if and only if $\Kk_2(2,A)$ is contained in the subgroup of $\St(2,A)$ 
		generated by $h_{12}(u)$, $h_{21}(u)$, $u\in \aa$ and $\Vv_2(A)=[\aa,\aa]$ $($the commutator subgroup of 
		$\aa$$)$. If $A$ is commutative, then $A$ is universal for $\GE_2$ if and only if $\Kk_2(2, A)$ is generated 
		by the Steinberg symbols.
	\end{thm}
	
	This result was known to experts of the field, but apparently no proof of it exists in the literature.
	Recently Hutchinson gave a proof of the second part of the above theorem (see \cite[App. A]{hutchinson2022}).
	
	\begin{prp}[Hutchinson]\label{hut}
		Let $A$ be a commutative ring. Then the naturl map $\St(2,A) \arr \Cc(A)$ given by $x_{12}(a)\mapsto \varepsilon(-a)\varepsilon(0)^3$ 
		and $x_{21}(a)\mapsto\varepsilon(0)^3\varepsilon(a)$ induces isomorphisms 
		\[
		\displaystyle \frac{\St(2,A)}{\Cc(2, A)}\simeq \Cc(A),\ \ \ \ \displaystyle \frac{\Kk_2(2,A)}{\Cc(2, A)}\simeq \Uu(A).
		\]
		In particular, $A$ is universal for $\GE_2$ if and only if $\Kk_2(2, A)$ is generated by Steinberg symbols.
	\end{prp}
	\begin{proof}
		See \cite[Theorem A.14, App. A]{hutchinson2022}.    
	\end{proof}
	
	\begin{exa}
		\par (i) (Morita) Let $A$ be a Dedekind domain, and $p \in A$ a nonzero prime element. Suppose that 
		\[
		\aa \arr (A/pA)^\times
		\]
		is surjective. If $\Kk_2(2,A)$ is generated by Steinberg symbols, then the same 
		is true for $\Kk_2(2,A\pth)$ \cite[Theorem 3.1]{morita2020}. In particular, if $A$ is universal for 
		$\GE_2$, then $A\pth$ is universal for $\GE_2$. For example since the natural map $\z^\times \arr (\z/p)^\times$
		is surjective for $p=2,3$, $\z\Big[\frac{1}{2}\Big]$ and $\z\Big[\frac{1}{3}\Big]$ are universal for $\GE_2$.
		\par (ii) (Morita) Let $p$ be a prime number. Then $\Kk_2(2, \z\pth)$ is generated by Steinberg symbols if and 
		only if $p=2,3$ \cite[Theorem 5.8]{morita2020}. Thus $\z\pth$ is universal for $\GE_2$ if and only if $p=2,3$. 
		See also \cite[Example 6.13, Lemma 6.15]{hutchinson2022}.
		
		\par (iii) Let $m=p_1\cdots p_k$ be an integer such that $p_i$ are primes and $p_1< \cdots<p_k$. Then it 
		follows from (i) that $\Kk_2(2, \z\Big[ \frac{1}{m} \Big])$ is generated by Steinberg symbols whenever 
		$(\z/p_i)^\times$ is generated by the residue classes $\{-1, p_1,\dots, p_{i-1}\}$ for all $i \leq k$. 
		Observe that the map $\z^\times  \arr (\z/p)^\times$ is surjective only for $p=2,3$. Thus $p_1$ in 
		above chain should be either $2$ or $3$. For example $\z\Big[\frac{1}{6}\Big]$, $\z\Big[\frac{1}{10}\Big]$, 
		$\z\Big[\frac{1}{15}\Big]$ and $\z\Big[\frac{1}{66}\Big]$ are universal for $\GE_2$.
		\par (iv) Let $k$ be a field. If $A$ is either $k[T]$ or $k[T,T^{-1}]$, then $\Kk_2(2,A)$ is generated by 
		Steinberg symbols \cite[Theorem 6]{Silvester1973}, \cite[Theorem 1]{alperin1979}.
	\end{exa}
	
	Let $a, b\in A$ be any two elements such that $1-ab\in \aa$. We define
	\[
	\lan a,b\ran_{ij}:=x_{ji}\bigg(\frac{-b}{1-ab}\bigg)x_{ij}(-a)x_{ji}(b)x_{ij}\bigg(\frac{a}{1-ab}\bigg)h_{ij}(1-ab)^{-1}.
	\]
	If $a$ and $b$ commute, then $\lan a,b\ran_{ij}\in \Kk_2(2, A)$. We call $\lan a,b\ran_{ij}$ a 
	{\it Dennis-Stein symbol}.
	If $u,v\in\aa$ commute, then 
	\[
	\displaystyle\{u,v\}_{ij}=\bigg\lan u, \frac{1-v}{u}\bigg\ran_{ij}=\bigg\lan \frac{1-u}{v},v\bigg\ran_{ij}.
	\]
	Hence over commutative rings, Dennis-Stein symbols generalize Steinberg symbols.
	
	\section{The abelianization of \texorpdfstring{$\Ee_2(A)$}{Lg}}
	
	The following theorem is the main result of this paper.
	
	\begin{thm}\label{H11}
		Let $A$ be a ring and let $M$ be the additive subgroup of $A$ generated by
		$axa-x$ and $\sum_{i=1}^l 3(b_i+1)(c_i+1)$, where  $x\in A$ and $a,b_i,c_i \in \aa$ provided that $\prod_{i=1}^l [b_i,c_i]=1$. 
		Then there is an exact sequence of abelian groups
		\[
		H_2(\Ee_2(A), \z) \arr \displaystyle \Uu(A)/[\Uu(A), \Cc(A)]
		\overset{\alpha}{\larr} A/M \overset{\beta}{\larr} \Ee_2(A)^\ab \arr 0,
		\]
		where $\alpha$ is induced by the map $\Cc(A) \arr A/M$, $\ep(x)\mapsto x-3$,
		and $\beta$ is induced by $y\mapsto E_{12}(y)$. Moreover, if $A$ is commutative,
		then this exact sequence is an exact sequence of $\z[\aa/(\aa)^2]$-modules.
	\end{thm}
	\begin{proof}
		From the Lyndon/Hochschild-Serre spectral sequence associated to the extension
		\[
		1 \arr \Uu(A) \arr \Cc(A) \arr \Ee_2(A) \arr 1,
		\]
		we obtain the five term exact sequence
		\[
		H_2(\Cc(A),\z) \!\arr\! H_2(\Ee_2(A),\z)\! \arr\! H_1(\Uu(A),\z)_{\Ee_2(A)} 
		\!\arr\! H_1(\Cc(A),\z) \!\arr\! H_1(\Ee_2(A),\z)\!\arr\! 0
		\]
		(see \cite[Corollary 6.4, Chp. VII]{brown1994}). For the middle term we have 
		\[
		H_1(\Uu(A),\z)_{\Ee_2(A)}\simeq \displaystyle \bigg(\frac{\Uu(A)}{[\Uu(A),\Uu(A)]}\bigg)_{\Ee_2(A)}\simeq 
		\displaystyle \frac{\Uu(A)}{[\Uu(A), \Cc(A)]}.
		\]
		We show that $H_1(\Cc(A),\z)\simeq A/M$. Consider the map
		\[
		\phi: \Cc(A) \arr A/M,  \ \ \ \ \ \prod\ep(a_i) \mapsto \sum (a_i-3).
		\]
		This map is well-defined. First note that in $A/M$ we have $a=a^{-1}$ and $12=0$. Thus
		\[
		\phi(h(a))=(-a-3)+(-a^{-1}-3)+(-a-3)=-3a-9=-3a+3=-3(a-1).
		\]
		Now we have
		\begin{align*}
			& \phi(\ep(x)\ep(0)\ep(y))=x-3+0-3+y-3=x+y-9=x+y+3,\\
			&\phi(h(-1)\ep(x+y))=-3(-1-1)+x+y-3=x+y+3,\\
			&\phi(\ep(x)h(a))=x-3-3(a-1)=x-3a, \\ 
			&\phi(h(a^{-1})\ep(axa))=-3(a^{-1}-1)+axa-3=-3(a-1)+x-3=x-3a.
		\end{align*}
		Moreover,
		\begin{align*}
			\phi(h(ab)h(a^{-1})h(b^{-1}))&=-3(ab-1)-3(a^{-1}-1)-3(b^{-1}-1)\\
			&=-3(ab-1)-3(a-1)-3(b-1)\\
			&=-3(ab+a+b+1)\\
			&=-3(a+1)(b+1).
		\end{align*}
		These show that the map $\phi$ is a well-defined homomorphism. Since $A/M$ is an abelian group,
		we have the homomorphism
		\[
		\bar{\phi}: \Cc(A)/[\Cc(A), \Cc(A)] \arr A/M, \ \ \ \ \ep(x) \mapsto x-3.
		\]
		Now we define
		\[
		\psi: A/M \arr \Cc(A)/[\Cc(A), \Cc(A)],  \ \ \ \ \ x\mapsto \ep(x)\ep(0)^{-1}.
		\]
		We show that this map is a well-defined homomorphism. Consider the items (i), (ii) and (iii)
		from the definition of $\Cc(A)$ (Section \ref{CA}). If in (i) we put $y=-x$, then
		we get 
		\[
		\ep(x)\ep(0)\ep(-x)=h(-1)\ep(0).
		\]
		Thus in $\Cc(A)/[\Cc(A), \Cc(A)]$, we have
		$h(-1)\ep(x)\ep(-x)=1$. From this we obtain
		\begin{align*}
			h(a)^2 &=h(-1)h(a)h(-a)\\
			&=h(-1)\ep(-a)\ep(-a^{-1})\ep(-a)\ep(a)\ep(a^{-1})\ep(a)\\
			&=\Big(h(-1)\ep(a)\ep(-a)\Big)^2 h(-1)\ep(a^{-1})\ep(-a^{-1})=1.
		\end{align*}
		Now we have
		\begin{align*}
			\psi(axa)&=\ep(axa)\ep(0)^{-1}=h(a)\ep(x)h(a)\ep(0)^{-1}
			=h(a)^2\ep(x)\ep(0)^{-1}=\ep(x)\ep(0)^{-1}=\psi(a).
		\end{align*}
		Moreover using (ii) for $x=0$ in $\Cc(A)/[\Cc(A), \Cc(A)]$, we have
		\[
		\ep(a)=h(a)\ep(a^{-1})h(a)=h(a)^2\ep(a^{-1})=\ep(a^{-1}).
		\]
		This implies that $h(-a)=\ep(a)\ep(a^{-1})\ep(a)=\ep(a)^3$ and hence
		\[
		h(a)=h(-1)\ep(a)^3=h(-1)\ep(a^{-1})^3.
		\]
		Furthermore by (i) we have $\ep(3x)=h(-1)\ep(x)^3$. Using this formula we obtain
		\begin{align*}
			\ep(3(a+1)(b+1)) &=\ep(0)\ep(ab)^3\ep(a)^3\ep(b)^3\ep(1)^3\\
			&=\ep(0)\ep(ab)^3\ep(a)^3\ep(b)^3h(-1)\\
			&=\ep(0)h(-1)\ep(ab)^3h(-1)\ep(a^{-1})^3h(-1)\ep(b^{-1})^3\\
			&=\ep(0)h(ab)h(a^{-1})h(b^{-1}).
		\end{align*}
		Thus
		\begin{align*}
			\psi(3(a+1)(b+1))&=\ep(3(a+1)(b+1))\ep(0)^{-1}
			=h(ab)h(a^{-1})h(b^{-1}).
		\end{align*}
		This shows that $\psi$ is well-defined. Since
		\[
		\psi(x+y)=\ep(x+y)\ep(0)^{-1}=h(-1)\ep(x)\ep(0)\ep(y)\ep(0)^{-1}=h(-1)\ep(x)\ep(y)=\psi(x)\psi(y),
		\]
		$\psi$ is a homomorphism of groups. Now it is easy to see that $\bar{\phi}$ and $\psi$ are 
		inverses of each other. Thus $\bar{\phi}$ is an isomorphism. This shows that the desired sequence is exact.
		
		Now let $A$ be commutative. We know that $\aa$ acts as follows on $\Ee_2(A)$:
		\[
		a.E(x):=\diag(a,1)E(x)\diag(a^{-1},1)=D(a)E(a^{-1}x).
		\]
		Now let us define the following action of $\aa$ on $\Cc(A)$:
		\[
		a.\varepsilon(x):=h(a)\varepsilon(a^{-1}x),
		\]
		\[
		a.(\varepsilon(x)\varepsilon(y)):=\Big(a.\varepsilon(x)\Big)\Big(a.\varepsilon(y)\Big).
		\]
		Observe that 
		\[
		a.(\varepsilon(x)\varepsilon(y))=\varepsilon(ax)\varepsilon(a^{-1}y).
		\]
		We show that this is in fact an action. Clearly $1.\varepsilon(x)=\varepsilon(x)$. Moreover, for any $a,b \in \aa$ we have
		\begin{align*}
			a.(b.\varepsilon(x)) & =a.\Big(h(b)\varepsilon(b^{-1}x)\Big)\\
			& =a.\Big(\varepsilon(-b)\varepsilon(-b^{-1})\varepsilon(-b)\varepsilon(b^{-1}x)\Big)\\
			& =\Big(a.\varepsilon(-b)\Big) \Big(a.\varepsilon(-b^{-1})\Big)\Big(a.\varepsilon(-b)\Big)\Big(a.\varepsilon(b^{-1}x)\Big)\\
			&= \varepsilon(-ab) \varepsilon(-(ab)^{-1})\varepsilon(-ab)\varepsilon((ab)^{-1}x)\\
			&=h(ab)\varepsilon((ab)^{-1}x)\\
			&=(ab). \varepsilon(x).
		\end{align*}
		This shows that the above action is well-defined.
		
		Clearly the natural map $\Cc(A) \arr \Ee_2(A)$ respects the above actions. Therefore the group $\aa$ acts naturally
		on the Lyndon/Hochschild-Serre spectral sequence of the extension $1 \arr \Uu(A) \arr \Cc(A) \arr \Ee_2(A) \arr 1$:
		\[
		E_{p,q}^2=H_p(\Ee_2(A),H_q(\Uu(A),\z))\Rightarrow H_{p+q}(\Cc(A),\z).
		\]
		This induces an action of $\aa$ on the five term exact sequence discussed in the beginning of this proof.
		
		Now we study the action of $\aa^2$ on the terms of this exact sequence. For $\Ee_2(A)$ we have
		\begin{align*}
			a^2.E(x) &:=\diag(a^2,1)E(x)\diag(a^{-2},1)\\
			& =D(a)(aI_2)E(x)(a^{-1}I_2) D(a)^{-1}\\
			& =D(a)E(x)D(a)^{-1}.
		\end{align*}
		Since $D(a)\in \Ee_2(A)$ and since the conjugation action induces a trivial action on homology groups 
		\cite[Proposition 6.2, Chap. II]{brown1994}, the action of $\aa^2$ on homology groups $H_k(\Ee_2(A),\z)$ 
		is trivial. The action of $\aa^2$ on $\Cc(A)$ also is induced by a conjugation:
		\begin{align*}
			a^2.\varepsilon(x)&=h(a^2)\varepsilon(a^{-2}x)\\
			&=h(a^2)h(a^{-1})\varepsilon(x)h(a^{-1})\\
			&=h(a)\varepsilon(x)h(a)^{-1}.
		\end{align*}
		This shows that the action of $\aa^2$ on homology groups $H_k(\Cc(A),\z)$ is trivial. For example through the 
		isomorphism $H_1(\Cc(A),\z)\simeq A/M$, on sees that $\aa$ acts on $A/M$ by the formula $a.\overline{x}=\overline{ax}$. Now we have
		\[
		a^2.\overline{x}=\overline{a^2x}=\overline{x(a^2-1)}+\overline{x}=\overline{x}.
		\]
		Finally if $\overline{X}\in \Uu(A)/[\Uu(A), \Cc(A)]$, then
		\[
		a^2.\overline{X}=\overline{h(a)Xh(a)^{-1}}=\overline{h(a)Xh(a)^{-1}X^{-1}}\ \overline{X}=\overline{X}.
		\]
		Thus $\aa^2$ acts trivially on the terms of the above exact sequence. This implies that the discussed sequence
		is an exact sequence of $\z[\aa/(\aa)^2]$-modules. This completes the proof of the theorem.
	\end{proof}
	
	\begin{rem}
		If $A$ is commutative, then by Proposition \ref{hut} we have the isomorphism
		\[
		\frac{\Uu(A)}{[\Uu(A),\Cc(A)]}\simeq \frac{\Kk_2(2,A)}{[\Kk_2(2,A), \St(2,A)]\Cc(2,A)}.
		\]
	\end{rem}
	
	\begin{cor}[Cohn \cite{cohn1968}]\label{universal1}
		Let $A$ be universal for $\GE_2$ and $M$ the additive subgroup of $A$ as described in the
		above theorem.
		Then $\Ee_2(A)^\ab\simeq A/M$.
	\end{cor}
	\begin{proof}
		Since $A$ is universal for $\GE_2$, $\Uu(A)=1$. Thus the claim follows from the above theorem.
	\end{proof}
	
	\begin{exa}\label{Z-n}
		(i) If $\aa=\{1,-1\}$, then $A/M=A/12\z$. Thus
		$\Ee_2(A)^\ab$ is a quotient of $A/12\z$. In particular,
		if $A$ is universal for $\GE_2$ and $\aa=\{1,-1\}$, then $\Ee_2(A)^\ab\simeq A/12\z$.
		Now by Example~\ref{example2}(vi),
		\[
		\Ee_2(\z[X_1, \dots, X_n])^\ab\simeq \z[X_1, \dots, X_n]/12\z \simeq\z/12 \oplus \bigoplus_{i\in \N} \z.
		\]
		\par (ii) Let $6\in \aa$. Then
		for any $x\in A$,
		\[
		x=3(2x2-x)-(3x3-x)\in M.
		\]
		This shows that $M=A$ and thus $\Ee_2(A)^\ab=1$.
	\end{exa}
	
	\begin{cor}\label{DS-cor}
		For any commutative ring $A$, we have the exact sequence
		\[
		\Kk_2(2, A)/\Cc(2, A) \arr A/M \arr \Ee_2(A)^\ab \arr 1.
		\]
		In particular, if $\Kk_2(2,A)$ is generated by Dennis-Stein symbols, then $\Ee_2(A)^\ab\simeq  A/N$, where $N$
		is the additive subgroup 
		\[
		N:=M+\lan de(d+e-3):d,e\in A \ \text{such that} \ 1-de\in \aa\ran.
		\]
	\end{cor}
	\begin{proof}
		The exact sequence follows from Theorem \ref{H11} and Proposition \ref{hut}.
		It is straightforward to check that through the composition
		\[
		\St(2, A) \arr \Cc(A) \arr A/M
		\]
		we have
		\[
		x_{12}(r) \mapsto -r,  \ \ \ x_{21}(r)\mapsto r,  \ \ \ h_{12}(a)\mapsto -3(a-1).
		\]
		Now if $\lan d,e\ran_{12}\in  \Kk_2(2,A)$ is a Dennis-Stein symbol, then under this composition we have
		\begin{align*}
			\lan d,e\ran_{12} &\mapsto -e(1-de)^{-1}+d+e-d(1-de)^{-1}+3(1-de-1)\\
			&= -e(1-de)+d+e-d(1-de)+3(1-de-1)\\
			&= +de(d+e-3).
		\end{align*}
		This completes the proof.
	\end{proof}
	
	\begin{cor}\label{H-1}
		Let $A$ be a ring such that $2\in \aa$. Then we have the exact sequence 
		\[
		H_2(\Ee_2(A),\z) \arr \displaystyle\frac{\Uu(A)}{[\Uu(A),\Cc(A)]} \arr 
		\displaystyle\frac{A}{\lan axa-x:a\in \aa, x\in A\ran} \arr \Ee_2(A)^\ab \arr 1.
		\]
		In particular, if $A$ is commutative and $2\in \aa$, then we have the exact sequence
		\[
		\Kk_2(2,A)/\Cc(2,A) \arr A/I\arr \Ee_2(A)^\ab \arr 1,
		\]
		where $I$ is the ideal generated by the elements $a^2-1$, $a\in\aa$.
	\end{cor}
	\begin{proof}
		Since $3=2^2-1$, for any $x\in A$ we have 
		\[
		3x=(2^2-1)x=2x2 -x \in \lan axa-x:a\in \aa, x\in A\ran.
		\]
		This implies that 
		\[
		M=\lan axa-x:a\in \aa, x\in A\ran.
		\]
		Thus the claim follows from Theorem~\ref{H11}.
	\end{proof}
	
	\section{The abelianization of \texorpdfstring{$\Ee_2(A)$}{Lg} for certain rings}
	
	Let $A$ be a local ring with maximal ideal $\mmm_A$. If $\char(A/\mmm_A)\neq 2, 3$, then
	$6\in \aa$. Thus by Example \ref{Z-n}(ii) we have $\Ee_2(A)^\ab=1$. The following proposition
	gives the precise structure of $\Ee_2(A)^\ab$ over commutative local rings.
	
	\begin{prp}\label{exa-2}
		Let $A$ be a commutative local ring with maximal ideal $\mmm_A$. Then
		\[
		\SL_2(A)^\ab=\Ee_2(A)^\ab\simeq
		\begin{cases}
			A/\mmm_A^2 &  \text{if $|A/\mmm_A|=2$}  \\
			A/\mmm_A &  \text{if $|A/\mmm_A|=3$.}  \\
			0 &  \text{if $|A/\mmm_A|\geq  4$}  \\
		\end{cases}
		\]
	\end{prp}
	\begin{proof} 
		If $|A/\mmm_A|\geq 4$, then there is  $a\in \aa$ such that $a^2-1\in \aa$. This implies that $A=M$ 
		and thus $\Ee_2(A)^\ab=1$. Let $A/\mmm_A\simeq\F_2$. Then $2\in \mmm_A, 3\in \aa$. Moreover, for 
		$a\in \aa$, $a\pm 1\in \mmm_A$. Thus $M\se \mmm_A^2$. Now if $a, b\in \mmm_A$, then  
		\[
		ab=3(a/3)b=3\Big[((a/3)-1)+1\Big]\Big[((b-1)+1\Big]\in M.
		\]
		This implies that $\mmm_A^2\se M$. Thus $M =\mmm_A^2$ and we have $\Ee_2(A)^\ab\simeq A/\mmm_A^2$.
		Finally let $A/\mmm_A\simeq\F_3$. If $a\in \mmm_A$, then $a-1, a-2\in \aa$. Since
		\[
		a=(a-2)^{-1}((a-1)^2-1)\in M,
		\]
		we have $\mmm_A\se M$. Clearly $M\se \mmm_A$. Therefore $\Ee_2(A)^\ab\simeq A/\mmm_A\simeq \z/3$.
	\end{proof}
	
	\begin{exa}\label{F2}
		Let $A$ be a commutative local ring such that $A/\mmm_A\simeq \F_2$. 
		Since $\mmm_A/\mmm_A^2$ is a $\F_2$-vector space, from the exact sequence
		\[
		0 \arr \mmm_A/\mmm_A^2 \arr A/\mmm_A^2 \arr \F_2 \arr 0
		\]
		it follows that
		\[
		|\Ee_2(A)^\ab|=1+\dim_{\F_2} (\mmm_A/\mmm_A^2).
		\]
		Let $A$ be a discrete valuation ring. Then $\mmm_A$ is generated by a prime element and 
		thus $\mmm_A/\mmm_A^2\simeq \F_2$. Therefore either $A/\mmm_A^2\simeq \z/2\times\z/2$,
		or $A/\mmm_A^2\simeq \z/4$. For example, if $A=\F_2[X]/\lan X^2\ran$, then  $\mmm_A^2=0$ 
		and thus $\Ee_2(A)^\ab\simeq A\simeq \z/2\times \z/2$. If $p$ is a prime and $A=\z_{(p)}$, 
		then $\Ee_2(A)^\ab\simeq \z/4$.
	\end{exa}
	
	\begin{exa}
		Let $A$ be a local ring not necessary commutative. Let $A':=Z(A)$ be the center of $A$. It is known that
		$A'$ is a local ring. If $|A'/\mmm_{A'}|\geq 4$, then as in above proposition we have $\Ee(A)^\ab=1$: Let
		$a\in {(A')}^\times$ such that $a^2-1\in {(A')}^\times$. If $y \in A$ and $y':=(a^2-1)^{-1}y$, then
		$y=ay'a-y' \in M$.
	\end{exa}
	
	Understanding the structure of the homology groups of $\SL_2(\z[\frac{1}{m}])$ have been topic of many articles, see 
	for example \cite{adem-naffah1998}, \cite{ww1998}, \cite{ae2014}. The following result completely calculates the first 
	integral homology of these groups.
	
	\begin{prp}
		Let $m$ be a square free natural number. Then
		\[
		\begin{array}{c}
			\SL_2(\z[\frac{1}{m}])^\ab=\Ee_2(\z[\frac{1}{m}])^\ab\simeq 
		\end{array}
		\begin{cases}
			0   & \text{if $2|m,\ 3|m$}  \\
			\z/3   & \text{if $2|m,\ 3\nmid m$}\\
			\z/4   & \text{if $2\nmid m,\ 3|m$}\\
			\z/12 &  \text{if $2\nmid  m,\ 3\nmid m$.}
		\end{cases}
		\]
	\end{prp}
	\begin{proof}
		Let $A_m:=\z[\frac{1}{m}]$. Note that $A_m^\times=\{\pm n^i: n\mid m, i\in \z\}$. Since $A_m$ is euclidean, 
		we have $\Ee_2(A_m)=\SL_2(A_m)$. 
		
		If $2|m$ and $3|m$, then $6\in A_m^\times$. Thus by Example \ref{Z-n}(ii), we have $\Ee_2(\z[\frac{1}{m}])^\ab=1$.
		
		Let $2|m$ and $3\nmid m$. Then $2\in A_m^\times$ and hence $3=2^2-1\in M$. This implies that $3A_m\se M$.
		On the other hand, for any $a,b\in A_m^\times$, clearly $3(a+1)(b+1)\in 3A_m$. Now consider
		$n^i\in A_m^\times$, $i\in \z$. Since $3\nmid n$, we have $3|n^2-1$. Therefore $3| (n^i)^2-1$ for any
		$i\in \z$. This implies that $M \se 3A_m$. Thus $3A_m=M$. Now it is easy to see that $A_m/M\simeq \z/3$.
		The inclusion $A_m\se \z_{(3)}$, induces the commutative diagram
		\[
		\begin{tikzcd}
			\z/3 \ar[r, "\simeq"]\ar[d, equal] & A_m/M \ar[d] \ar[r, two heads]& \Ee_2(A_m)^\ab\ar[d]\\
			\z/3 \ar[r, "\simeq"] & \z_{(3)}/M \ar[r, "\simeq "]& \Ee_2(\z_{(3)})^\ab.
		\end{tikzcd}
		\]
		By Proposition \ref{exa-2}, the bottom maps are isomorphisms. Thus the upper right map is injective.
		This proves that $\Ee_2(A_m)^\ab\simeq \z/3$.
		
		Let $2\nmid m$ and $3\mid m$. Note that $3\in A_m^\times$. We show that $M=4A_m$. First note that
		\[
		4=12-8=3(1+1)(1+1) -(3^2-1)\in M.
		\]
		Since for any $i\geq 0$,
		\[
		4m^i=12m^i -m^i(3^2-1),\ \ \ \ 
		4/m^i=4m^i-\frac{4}{m^i}((m^i)^2-1)=12 m'-\frac{4}{m^i}((m^i)^2-1),
		\]
		where $m'\in\z$, we have $4A_m\in M$. On the other hand, let $n\mid m$. Since $n$ is odd, $2$ divides 
		$n-1$ and $n+1$. Thus for any $i,j\in\z$, $4\mid (n^i)^2-1$ and $4\mid 3(\pm n^i+1)(\pm n^j+1)$. These 
		facts imply that $M\se 4A_m$. Therefore $M=4A_m$. Now one easily verifies that
		\[
		A_m/M=A_m/4A_m\simeq \z/4.
		\]
		Since $2\nmid m$, $A_m \se \z_{(2)}$. Now with an argument similar to the previous case, using 
		Proposition \ref{exa-2} and Example \ref{F2},  one can show that $\Ee_2(A_m)^\ab\simeq \z/4$.
		Finally let $2\nmid m$ and $3\nmid m$. As previous cases we can show that $M=12A_m$. Thus
		\[
		A_m/M=A_m/12A_m\simeq \z/12.
		\]
		Since $2\nmid m$ and $3\nmid m$, we have $A_m\se \z_{(2)}$ and $A_m\se \z_{(3)}$. Now from the 
		commutative diagram
		\[
		\begin{tikzcd}
			\z/12 \ar[r, "\simeq"]\ar[d, two heads] & A_m/M \ar[d, two heads] \ar[r, two heads]& \Ee_2(A_m)^\ab\ar[d]\\
			\z/4 \ar[r, "\simeq"] & \z_{(2)}/M \ar[r, "\simeq "]& \Ee_2(\z_{(2)})^\ab
		\end{tikzcd}
		\ \ \ \ \ \ 
		\begin{tikzcd}
			\z/12 \ar[r, "\simeq"]\ar[d, two heads] & A_m/M \ar[d, two heads] \ar[r, two heads]& \Ee_2(A_m)^\ab\ar[d]\\
			\z/3 \ar[r, "\simeq"] & \z_{(3)}/M \ar[r, "\simeq "]& \Ee_2(\z_{(3)})^\ab
		\end{tikzcd}
		\]
		it follows that the composition 
		\[
		\z/12 \overset{\simeq}{\larr} A_m/M \two \Ee_2(A_m)^\ab
		\]
		is injective. This completes the proof of the proposition.
	\end{proof}
	
	\begin{cor}
		Let $m$ be a square free natural number. Then
		\[
		\begin{array}{c}
			\PSL_2(\z[\frac{1}{m}])^\ab\simeq 
			\begin{cases}
				0 & \text{if $2\mid m$, $3\mid m$}\\
				\z/3 & \text{if $2\mid m$, $3\nmid m$.}\\
				\z/2 & \text{if $2\nmid m$, $3\mid m$}\\
				\z/6 & \text{if $2\nmid m$, $3\nmid m$}
			\end{cases}
		\end{array}
		\]
	\end{cor}
	\begin{proof}
		Let $A_m=\z[\frac{1}{m}]$. From the extension 
		\[
		1 \arr \mu_2(A_m) \arr \SL_2(A_m) \arr  \PSL_2(A_m) \arr 1,
		\]
		we obtain the exact sequence
		\[
		\mu_2(A_m) \arr \SL_2(A_m)^\ab \arr  \PSL_2(A_m)^\ab \arr 1.
		\]
		Now the claim follows from the above proposition and the fact that 
		\[
		\PSL_2(A_m)^\ab\simeq H_1(\PSL_2(A_m),\z)\simeq (\z/2)^\alpha \oplus (\z/3)^\beta
		\]
		for some $\alpha$ and $\beta$ by \cite[Corollary 4.4]{ww1998}.
	\end{proof}

	Let $\OO_d$ be the ring of algebraic integers of the quadratic field $\q(\sqrt{d})$, 
	where $d$ is a square-free integer. It is known that
	\[
	\OO_d=\begin{cases}
		\z[\sqrt{d}] & \text{if $d\equiv 2,3\pmod 4$}\\
		\z[\frac{1+\sqrt{d}}{2}] & \text{if $d\equiv 1\pmod 4$.}\\
	\end{cases}
	\]
	
	If $d<0$, then $\OO_d^\times$ has at most 6 elements.  In fact 
	\[
	\OO_{-1}^\times=\{\pm 1, \pm i\}, \ \ \ \ 
	\OO_{-3}^\times=\{\pm 1, \pm \omega , \pm (\omega - 1)\},
	\]
	where $i=\sqrt{-1}$, $\omega=\displaystyle\frac{1+\sqrt{-3}}{2}$
	and for $d\neq -1, -3$, 
	\[
	\OO_d^\times=\{\pm 1\}.
	\]
	It is known that $\OO_d$ is norm-euclidean if and only if $d = -1,-2,-3,-7,-11$ if and only if 
	\[
	\Ee_2(\OO_d)=\SL_2(\OO_d)
	\]
	(see \cite[Theorem 5.1]{ELS1992} and \cite[Theorem 6.1]{cohn1966}). Observe that for $d<0$, 
	$\OO_d$ is universal for $\GE_2$ if and only if $d\neq-2,-7,-11$ (see Example \ref{example2}(v) 
	and \cite[Remarks, pages 162-163]{cohn1968}). It is easy to see that
	\[
	M=\begin{cases}
		2\OO_{-1} & \text{if $d=-1$}\\
		(2\omega-1)\OO_{-3} & \text{if $d=-3.$}\\
		12\z & \text{otherwise}
	\end{cases}
	\]
	It follows from this that
	\[
	\OO_d/M\simeq\begin{cases}
		\z/2 \times\z/2 & \text{if $d=-1$}\\
		\z/3 & \text{if $d=-3.$}\\
		\z/12\times \z  & \text{otherwise}
	\end{cases}
	\]
	\begin{exa}
		If $d=-2,-7,-11$, then we have a surjective map
		\begin{equation}\label{two}
			\z/12 \times \z\simeq \OO_d/M \two \Ee_2(\OO_d)^\ab.
		\end{equation}
		\par (i) For $d=-2$, we have $1-(-\sqrt{-2})(\sqrt{-2})\in \OO_{-2}^\times$. Thus 
		$\lan -\sqrt{-2},\sqrt{-2} \ran_{12} \in \Kk_2(2, \OO_{-2})$ is a Dennis-Stein symbol. Hence 
		under the map (\ref{two}) the element
		\[
		-6=(-\sqrt{-2})(\sqrt{-2})\bigg(-\sqrt{-2}+\sqrt{-2}-3\bigg)=2(-3) \in \OO_{-2}/M
		\]
		maps to zero (see Corollary \ref{DS-cor}). Thus we have a surjective map $\z/6 \times \z \arr \Ee_2(\OO_{-2})^\ab$. 
		\par (ii) For $d=-7$, we have 
		$\displaystyle 1-\Big(\frac{1+\sqrt{-7}}{2}\Big)\Big(\frac{1-\sqrt{-7}}{2}\Big)\in \OO_{-7}^\times$. 
		Thus 
		\[
		\displaystyle\Big\lan \frac{1+\sqrt{-7}}{2}, \frac{1-\sqrt{-7}}{2} \Big\ran_{12} \in \Kk_2(2, \OO_{-7})
		\]
		is a Dennis-Stein symbol. Again under the map (\ref{two}) the element
		\[
		4=(\frac{1+\sqrt{-7}}{2})(\frac{1-\sqrt{-7}}{2})\bigg(\frac{1+\sqrt{-7}}{2}+\frac{1-\sqrt{-7}}{2}-3\bigg) \in\OO_{-7}/M
		\]
		maps to zero. Thus we have a surjective map  $\z/4 \times \z \arr \Ee_2(\OO_{-7})^\ab$. 
		\par (iii) Let $d=-11$. It is easy to see that in $\Kk_2(2, \OO_{-11})$ any Dennis-Stein symbol is a Steinberg symbol.
		In fact, there are no $x,y \in \OO_{-11}$, such that $1-xy=-1$ or equivalently $xy=2$. But if 
		$\displaystyle x=\frac{1+\sqrt{-11}}{2}$, then $x\overline{x}=3$. Now we have 
		\[
		(E(x)E(\overline{x}))^3=-I_2.
		\]
		From this we obtain the element
		\[
		\Theta:=h(-1)(\varepsilon(x)\varepsilon(\overline{x}))^3\in \Uu(\OO_{-11})\simeq \frac{\Kk_2(2, \OO_{-11})}{\Cc(2, \OO_{-11})}.
		\]
		Now it is straightforward to see that under the map 
		$\displaystyle\frac{\Kk_2(2, \OO_{-11})}{\Cc(2, \OO_{-11})} \arr \OO_{-11}/M$, we have
		\[
		\Theta\mapsto -3(-1-1)+3(x-3+ \overline{x}-3)=6-15=-9=3.
		\]
		Thus we have a surjective map $\z/3 \times \z \arr \Ee_2(\OO_{-11})^\ab$.
	\end{exa}
	The following theorem is due to P.M. Cohn.
	
	\begin{prp}[Cohn \cite{cohn1966},\cite{cohn1968}]\label{Od-}
		Let $d<0$ be a square free integer. Then 
		\[
		\Ee_2(\OO_d)^\ab \simeq 
		\begin{cases}
			\z/2 \times \z/2 & \text{if $d=-1$}\\
			\z/6 \times \z & \text{if $d=-2$}\\
			\z/3 & \text{if $d=-3$}\\
			\z/4 \times \z & \text{if $d=-7$}\\
			\z/3 \times \z & \text{if $d=-11$}\\
			\z/12 \times \z & \text{otherwise}
	\end{cases}.
	\]
\end{prp}
\begin{proof}
	We have seen that $\OO_d$ is universal for $\GE_2$ if and only if $d\neq -2,-7,-11$ (see Example \ref{example2}
	and \cite[Remarks, pages 162-163]{cohn1968}). Thus by Corollary \ref{universal1}, 
	\[
	\Ee_2(\OO_d)^\ab\simeq \OO_d/M.
	\]
	This proves the claim for $d\neq -2,-7,-11$. For the case $d= -2,-7,-11$, see \cite[page 162]{cohn1968}.
\end{proof}

Now let $d > 0$. Then $\OO_d^\times$ has infinitely many units. In fact $\OO_d^\times=\{\pm u^n:n\in \z\}$, 
where $u$ is a particular unit called a {\it fundamental unit}. For a fundamental unit $u$, there are three 
other fundamental units:  $\overline{u}$, $-u$ and $-\overline {u}$. In fact the one of these four elements 
which is greater than $1$ is called the ``fundamental unit". Observe that $\Ee_2(\OO_d)=\SL_2(\OO_d)$ 
\cite[page 321, Theorem]{vas1972}.

\begin{prp}
	Let $d>0$ be a square free integer. Let $u$ be the fundamental unit of $\OO_d^\times$. If
	$d\equiv 2, 3 \pmod 4$ and $u=a+b\sqrt{d}$, then
	\[
	\OO_d/M\simeq \begin{cases}
		\displaystyle\frac{\z}{2\gcd(bd, 3a+3, 6)} \times \frac{\z}{2b} & \text{if $N(u)=1$}\\
		\\
		\displaystyle\frac{\z}{2\gcd(a,3)} \times \frac{\z}{2\gcd(a, 3b)} & \text{if $N(u)=-1$}
	\end{cases}
	\]
	and if $d\equiv 1 \pmod 4$ and $u=a+b(1+\sqrt{d})/2$, then
	\[
	\OO_d/M\simeq 
	\begin{cases}
		\displaystyle\frac{\z}{\gcd(b, 6(a-1), 12)} \times \frac{\z}{\gcd(bd, 12(a-1)+6b, 24)}  & \text{if $N(u)=1$}\\
		\\
		\displaystyle\frac{\z}{\gcd(2a+b,6(a-1), 12)} \times \frac{\z}{\gcd(2a+b,  6b)} & \text{if $N(u)=-1$}.
	\end{cases}
	\]
	In particular, $\Ee_2(\OO_d)^\ab$ is a finite group.
\end{prp}
\begin{proof}
	Consider the isomorphism 
	\[
	\OO_d/M\simeq \Big(\OO_d /\lan u^2-1\ran\Big) /\Big(M/\lan u^2-1\ran\big).
	\]
	Then
	\[
	\lan u^2-1\ran=\lan \overline{u}(u^2-1)\ran=
	\begin{cases}
		\lan u- \overline{u} \ran& \text{if $N(u)=1$}\\
		\lan u + \overline{u} \ran& \text{if $N(u)=-1.$}
	\end{cases}
	\]
	First let $d\equiv 2, 3 \pmod 4$.  Then $\OO_d=\z[\sqrt{d}]$ and
	$u-\overline{u}=2b\sqrt{d}$ when $N(u)=1$ and $u+ \overline{u}=2a$ when $N(u)=-1$. Hence
	\[
	\OO_d /\lan u^2-1\ran=\begin{cases}
		\OO_d/\lan 2b\sqrt{d}\ran & \text{if $N(u)=1$}\\
		\OO_d/\lan 2a\ran & \text{if $N(u)=-1$}
	\end{cases}
	\simeq
	\begin{cases}
		\z/2bd \times \z/2b & \text{if $N(u)=1$}\\
		\z/2a \times \z/2a & \text{if $N(u)=-1.$}
	\end{cases}
	\]
	Since $\overline{u}^2=1$ in $\OO_d /\lan u^2-1\ran$, we have
	\[
	M/\lan u^2-1\ran=\lan \overline{6(u+1)}, \overline{12}\ran=
	\begin{cases}
		\lan \overline{6(a+1)},\overline{12}\ran & \text{if $N(u)=1$}\\
		\lan \overline{6(1+b\sqrt{d})}, \overline{12}\ran & \text{if $N(u)=-1$}.
	\end{cases}
	\]
	Thus
	\[
	\OO_d/M\simeq \begin{cases}
		\z/\gcd(2bd, 6(a+), 12) \times \z/2b & \text{if $N(u)=1$}\\
		\z/\gcd(2a,6)\ran \times \z/ 2\gcd(2a, 6b) & \text{if $N(u)=-1$}.
	\end{cases}
	\]
	Now let $d\equiv 1 \pmod 4$. Then $\OO_d=\z[\omega]$, where $\omega=(1+\sqrt{d})/2$.
	If $u=a+b\omega\in \OO_d$, then $\overline{u}=a+b\overline{\omega}$, where
	$\overline{\omega}=(1-\sqrt{d})/2$. Note that $u-\overline{u}=b\sqrt{d}$ if $N(u)=1$ 
	and $u+ \overline{u}=2a+b$ if $N(u)=-1$. Thus we have
	\begin{align*}
		\OO_d /\lan u^2-1\ran &=\begin{cases}
			\OO_d/\lan b\sqrt{d}\ran & \text{if $N(u)=1$}\\
			\OO_d/\lan 2a+b\ran & \text{if $N(u)=-1$}
		\end{cases}\\
		& \simeq
		\begin{cases}
			(\z \times \z)/\lan b(-1,2),b((d-1)/2,1)\ran  & \text{if $N(u)=1$}\\
			\z/(2a+b) \times \z/(2a+b) & \text{if $N(u)=-1$}
		\end{cases}\\
		& \simeq
		\begin{cases}
			\z/b \times \z/bd  & \text{if $N(u)=1$}\\
			\z/(2a+b) \times \z/(2a+b) & \text{if $N(u)=-1$},
		\end{cases}
	\end{align*}
	where the isomorphism 
	\[
	(\z \times \z)/\lan b(-1,2),b((d-1)/2,1)\ran \arr \z/b \times \z/bd
	\]
	is given by $\overline{(r,s)}\mapsto (\overline{r+s}, \overline{2r+s})$.
	Moreover, note that
	\begin{align*}
		M/\lan u^2-1\ran &=\lan \overline{6(u-1)}, \overline{12}\ran\\
		& =\lan \overline{6(a-1)+6b\omega}, \overline{12}\ran \\
		&\simeq\lan (\overline{6(a-1)},\overline{6b}),\overline{(12,0)}\ran\\
		&\simeq
		\begin{cases}
			\lan (\overline{6(a-1)},12\overline{(a-1)}+\overline{6b}),\overline{(12,24)}\ran & \text{if $N(u)=1$}\\
			\lan (\overline{6(a-1)},\overline{6b}),\overline{(12,0)}\ran & \text{if $N(u)=-1$}.
		\end{cases}
	\end{align*}
	Thus we obtain the desired isomorphism.
\end{proof}



\end{document}